\documentclass[11pt]{article}

\usepackage{cite}

\usepackage{graphicx}
\usepackage{amsfonts}
\usepackage{amssymb}
\usepackage{amsmath}
\usepackage{amsthm}

\def\N{\mathbb N}

\newtheorem{definition}{Definition}
\newtheorem{proposition}{Proposition}
\newtheorem{theorem}{Theorem}

\newtheorem{lemma}{Lemma}
\newtheorem{corollary}{Corollary}

\title{The problem of combinatorial encoding of a continuous dynamics and the notion of transfer of paths in graphs}

\author{A.~M.~Vershik\thanks{St.Petersburg Department
 of Steklov Institute of Mathematics and
 St.Petersburg State University,
St.Petersburg, Russia,
Institute for Information Transmission Problems,
Moscow, Russia. E-mail: {\tt avershik@pdmi.ras.ru}. Partially supported by the RFBR grant 17-01-00433.}}

\date{October 8, 2019}

\begin{document}
\maketitle

\begin{abstract}
We introduce the notion of combinatorial encoding of continuous dynamical systems and suggest the
first examples, which are the most interesting and important, namely, the combinatorial encoding of a Bernoulli process with continuous state space, e.g., a sequence of i.i.d. random variables with values in the interval with the Lebesgue measure (or a Lebesgue space).

The main idea is to associate with a random object (a trajectory of the random process) a path in an $\N$-graded graph and parametrize it with the vertices of the graph that belong to this path. This correspondence (encoding) is based on the definition of a decreasing sequence of cylinder partitions, and the first problem is to verify whether or not the given combinatorial encoding has the property of distinguishability, which means that our encoding is an isomorphism, or, equivalently, the limit of the increasing sequence of finite partitions is the partition into singletons $\bmod\,0$. This is a generalization of the problem of generators in ergodic theory.

The existence of a suitable $\N$-graded graph  is equivalent to the so-called standardness of the orbit partition in the sense of the theory of filtrations in measure spaces.

In the last section, we define  the notion of a so-called transfer, a transformation of paths in a graded graph, as a generalization of the shift in stationary dynamics.

\end{abstract}

\section{General problems}

\subsection{Combinatorial encoding}

We will consider $\N$-graded locally finite graphs $\Gamma$ (Bratteli diagrams) and the spaces of their
maximal infinite paths~$T(\Gamma)$; clearly, $T(\Gamma)$ endowed with the natural topology is a compact space. We will also consider a filtration $\{T_n\}$, which is a sequence of quotient spaces of $T(\Gamma)$. Namely, $T_n$ is the set of all paths that start at vertices of level $n$, $n=1,2, \dots$; here $T_1=T(\Gamma)$.
Thus, $T_n$ is the result of forgetting the first $n$ vertices of each path, or the quotient space of $T(\Gamma)$ with respect to the $n$-equivalence of paths (two paths are $n$-equivalent if they coincide starting from level~$n+1$). This means that we have a decreasing sequence of partitions $\xi_n$, $n=1,2 \dots$, and $T_n=T(\Gamma)/{\xi_n}$.

So, we have a chain of quotient spaces:
$$ T(\Gamma)=T_1\longrightarrow T_2\longrightarrow\dots. $$
This is a canonical filtration on the path space $T(\Gamma)$ of an $\N$-graded graph~$\Gamma$.
If we define a Borel measure on the space  $T(\Gamma)$, then we have a filtration on the resulting measure space. In particular, if the measure is  a central measure (see \cite{UMN}) on the path space, then $\{\xi\}_n$ is a semi-homogeneous filtration.

On the other hand, consider an infinite product of measure spaces (e.g., the unit intervals with the Lebesgue measure) in the ordinary sense:
 $$(I^{\infty}, m^{\infty})=\prod_{n=1}^{\infty}(I^n,m^n)$$
 (a Bernoulli space).

Assume that we have a hyperfinite equivalence relation $\tau$ on the space $(I^{\infty}, m^{\infty})$,
which is, by definition, the limit of a  sequence of equivalence relations $\tau_n$
where $\tau_n$ is an $n$-cylinder equivalence relation on $(I^{\infty}, m^{\infty})$ with finite blocks.

An important example appears if we have an action of a countable locally finite group
$G=\bigcup_n G_n$ and  $\tau$ is the orbit partition of $G$ while  $\tau_n$ is the orbit partition of the finite group $G_n$.

\begin{definition}
A homogeneous hyperfinite equivalence relation $\tau=\lim_n \tau_n$ is said to be \emph{standard} if there exists an increasing sequence of finite partitions $\{\eta_n\}$ such that

{\rm(i)} $\eta_n$ is an independent complement to $\tau_n$
and

{\rm (ii)} $\bigvee_n \eta_n=\epsilon$.

\noindent Here $\epsilon$ stands for the partition of the space into
singletons  $\bmod\,0$; in other words, the sequence of partitions $\{\eta_n\}_{n=1,2 \dots }$ separates almost all pairs of points.
\end{definition}

Both conditions are essential: there exists an increasing sequence that satisfies condition~(i) but does not satisfy condition~(ii).

For our purposes, it is convenient to define an increasing sequence of finite cylinder measurable partitions $\{\eta_n\}_{n=1,2 \dots }$ in the space $I^{\infty}$ (this means that $\eta_n$ is a  subpartition of the partition according to the first $n$ coordinates of the infinite product $I^{\infty}$ with the product measure $m^{\infty}$).

As a result, we have a map
$$
L:I^{\infty}\longrightarrow T(\Gamma)
$$
which sends the sequence $\{\eta_n\}_n$ to a canonical sequence on the path space of an $\N$-graded graph $\Gamma$, which is defined automatically if we assume that the $L$-image of the partition $\eta_n$ is the partition of the set of paths in $\Gamma$ according to the vertices of the $n$th level $\Gamma_n$ of the graph $\Gamma$, $n=1,2, \dots$; the edges of~$\Gamma$ are defined according to the structure of the partition $\eta_n$ with respect to the coordinates of the space $I^{\infty}$.
A concrete example of all these steps will be given in the next section.

In order to check condition~(ii), it suffices to prove that the map $L$ is an isomorphism
of measure spaces between $(I^{\infty}, m^{\infty})$ and $(T(\Gamma),\mu)$ where $\mu$ is
the central measure that coincides with the $L$-image of the measure~$m^{\infty}$. This is the main point of our considerations.

If both conditions are satisfied, then we say that the increasing sequence $\{\eta_n\}_n$ of finite partitions of $(I^{\infty},m^{\infty})$  is a \emph{combinatorial encoding}
of the equivalence relation $\tau$ on $(I^{\infty},m^{\infty})$ (or of the Bernoulli scheme), and the
distinguishabilitity problem for the sequence $\{\eta_n\}_n$ has a positive solution.

\subsection{Transfer}
 Our method of constructing the map $L$ in the case where $I^{\infty}=\prod_n I$ (for example, $I=[0,1]$)  allows us to define the transformation
 $LSL^{-1}$ of $T(\Gamma)$ as the image of the one-sided shift $S$ on $I^{\infty}$.  A~well-known example is the Sch\"utzenberger transformation on the Young graph~$\textbf{Y}$, which is the image of the one-sided shift under the left ($Q$) part of the RSK correspondence on $T(\textbf{Y})$, see~\cite{KV}. In \cite{RS, Sn},
it is proved that this map is an isomorphism. We will return to this subject in a separate paper. At the end of this paper, we give a direct definition of a \emph{transfer}, which is an endomorphism (shift) on the space $T(\Gamma)$,
 in graph terms for all graphs with the following property: each 2-interval of the graph has at most 2 intermediate vertices. A general example of such a graph is
 the Hasse diagram of a distributive lattice.

\subsection{Several remarks}

1. The classical encoding process corresponds to the case where $\Gamma$ is a stationary graph
(all levels, as well as all sets of edges between adjacent levels, are isomorphic).

\smallskip
2. The standard filtrations do not exhaust the class of all filtrations, and standard actions of locally finite groups are not even typical. In the generic case of a hyperfinite equivalence relation, there is no partition that satisfies both conditions (i)~and~(ii). For the theory of general
filtrations and the related literature, see~\cite{UMN}.

\smallskip
3. The list of $\N$-graded graphs that can appear in this context is a very interesting class of graphs. To my knowledge, there is no exact description of this class.

\smallskip
4. Let us mention some important examples.

\begin{enumerate}
\item[4A.] The action of the group ${\mathfrak S}_{\infty}$ on $([0,1]^{\infty},m^{\infty})$ by permutations (see the next section).

\item[4B.] The action of the group $U(\infty)$ on the space of infinite Hermitian matrices with Gaussian measure
(GUE, GOE); it will be considered in a paper under preparation by the author and F.~Petrov.

\item[4C.] An equivalence relation coming from the RSK algorithm (a paper under preparation)

\item[4D.] Cases  that should be considered: the action of ${\mathfrak S}_{\infty}$ on the space of  infinite $\{0,1\}$-matrices by conjugation, the actions of $U(\infty)$ and $O(\infty)$ on the spaces of forms in $m$ variables, etc.
\end{enumerate}

\smallskip
5. 
Let $Z$ be a countable set, $I=[0,1]$, $m$ be the Lebesgue measure on $I$, and
$X=I^Z$, $\mu=m^Z$.
Assume that we have a family of finite measurable partitions $\{\xi_F\}$
where $F$ is a finite subset of $Z$ and  $\xi_F$ is a finite partition of
$I^F$ which can be regarded as a cylinder partition of $X$.
A family $\{\xi_F\}$ is called distinguishable if the product of
the partitons  $\xi_F$ over all $F$
 is the partition of $X$ into singletons $\bmod \,0$ with respect to the
measure $\mu$.

\smallskip
\noindent{\bf Conjecture.} For all distinguishable
families of partitions $\{\xi_F\}$,
$$
\lim_n H(\xi_{F_n)}/|F_n|=\infty
$$
where $H(\cdot)$ is the entropy of
a partition, $|\cdot|$ is the number of elements, and the limit is
over a sequence of
finite subsets of $Z$ such that the union of $F_n$ is $Z$.
\smallskip

For one special case, this conjecture was proved by the student
G.~Veprev in the paper published in this volume.


\section{Example: the Weyl encoding of a Bernoulli process and the positive solution to the distinguishability problem}

Our first example is as follows. Let $\eta_n$ be a  \emph{cylinder partition} with the partition of the $n$-dimensional cube $I^n$ into open \emph{Weyl simplices} for the base. The  last sentence means the following:
we regard the space ${\Bbb R}^n$ as a Cartan subalgebra of the $A_n$ Lie algebra with a fixed root system. By an ``open Weyl simplex'' we mean
the intersection of an open Weyl chamber in the ordinary sense with the unit cube $[0,1]^n.$ We also assume that the correspondence between the Weyl chambers and the elements of the Weyl group (which is the symmetric group ${\mathfrak S}_n$) is also fixed. Although the language of the theory of Lie algebras is not necessary for describing our example, we nevertheless use it in order to generalize the construction  in future to the case of other series of simple Lie algebras.

 We will consider the subset of $I^n$ of full (Lebesgue) measure $m^n$  that consists of the vectors with pairwise distinct coordinates.

Each Weyl simplex is identified with some permutation of the set $\textbf{n}=\{1, 2, \dots, n\}$, an element of the Weyl group.

Recall that the fact that the partitions $\eta_n$ are \emph{finite} means that the number of blocks of positive measure in~$\eta_n$ is finite; in our
case, it is equal to $\operatorname{ord}({\mathfrak S}_n) = n!$.

The fact that the sequence $\{\eta_n\}_n$ is  \emph{increasing} means that for every~$n$ each element $C$ of the partition
$\eta_{n+1}$
(regarded as a subset of $I^n$) is a subset of an element $D$ of the partition $\eta_n$, and each element $D$ of the partition  $\eta_n$
(regarded as subset of $I^n$) is exactly the disjoint union of
all elements of $\eta_{n+1}$ that belong to it. The union of all elements of all partitions $\eta_n$, $n=1,2 \dots$ (the set of all Weyl simplices), generates a \emph{tree} which we denote by $W$ and call the permutation  (or factorial) tree.

The $n$th level of the tree $W$ has $n!$ vertices  which correspond to the elements of the partition $\eta_n$, and each vertex can be identified with a permutation from ${\mathfrak S}_n$.

The edges of the tree join each element  $C$ of $\eta_{n+1}$ with the element~$D$ of~$\eta_n$ that contains $C$; see Fig.~1.
We define coordinates for Weyl simplices and the corresponding permutations as follows.
Let $x^n=(x_1,x_2,\dots, x_n)$ be a point of $I^n$; then the coordinates
$$(k_1,k_2,\dots, k_n)$$
corresponding to the simplex $\sigma_{x^n}$
are given by the formula
    $$k_i=\#\{s \in {\bf n}: x_s<x_i\}, \quad i=1,2, \dots, n.$$
In particular, the coordinates $(k_1,k_2,\dots, k_n)$ do not depend on the choice of a point of the simplex.

\begin{figure}
\includegraphics[scale=0.7]{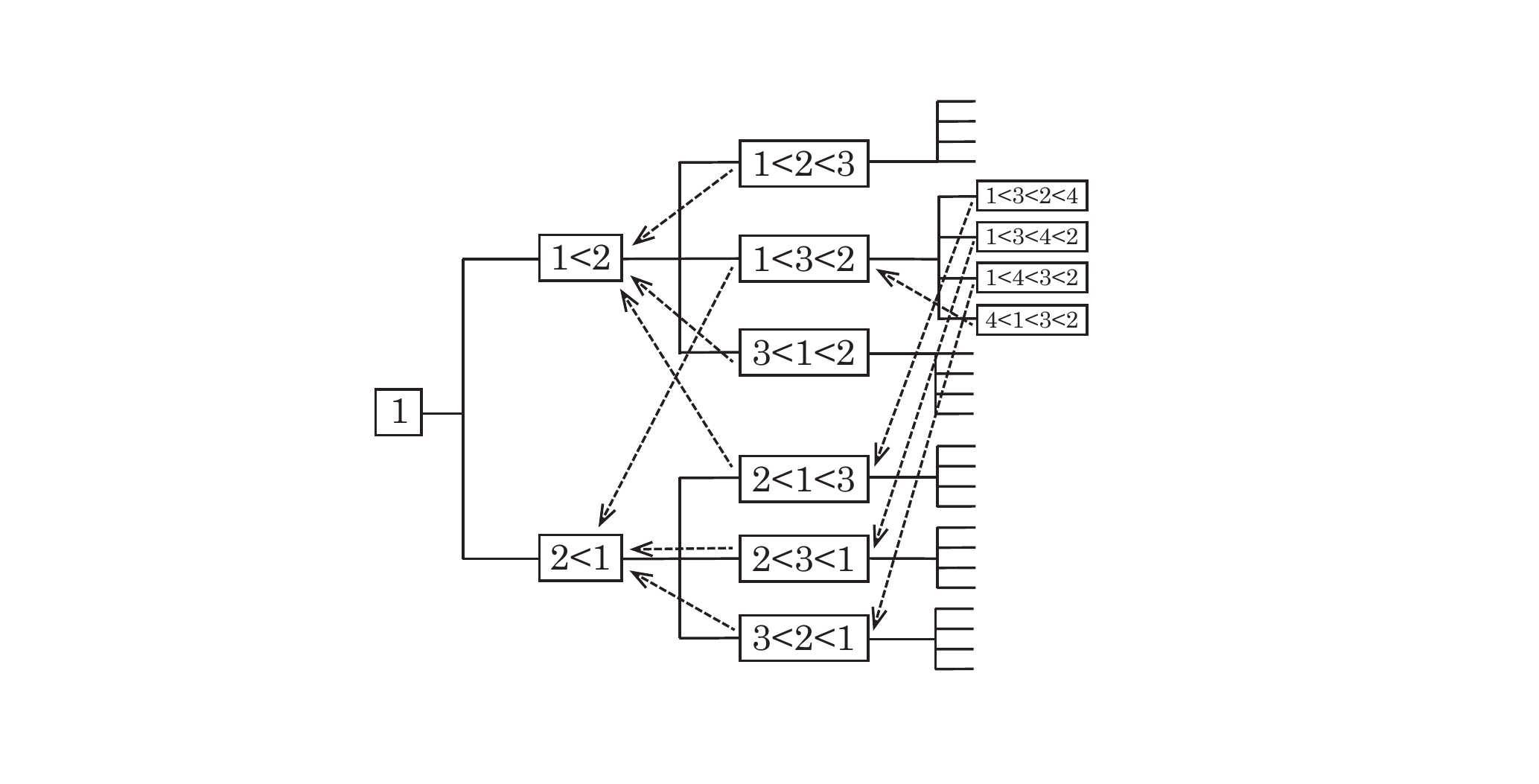}
\caption{A permutation tree with a translation.} \label{fig:Tree}
\end{figure}

The \emph{invariance} of the sequence of partitions $\{\eta_n\}_n$ with respect to the one-sided shift $S$ of the cube $I^{\infty}$ means that the images
$Sx, Sx'$ of almost all pairs of points $x, x'$ of $I^{\infty}$ that belong to the same element of the partition $\eta_n$ (or to the same Weyl simplex)
will belong to the same element of the partition $\eta_{n-1}$ (or to the same Weyl simplex) for $n = 2,3 \dots$.

This means that we have constructed another system of edges on the tree~$W$: each element of the partition $\eta_n$ that corresponds to a vertex of level $n$ has an edge which goes to the vertex of the previous level~$n-1$ that corresponds to some element of the partition $\eta_{n-1}$. Such edges between permutations will be called \emph{translations}.

\begin{proposition}
The permutation $(r_1,r_2,\dots,r_{n-1})$ corresponding to the shifted simplex
$\Sigma_{Sx^n}$ is given by the formula
      $$
 r_i=\begin{cases}
 k_{i+1} &\text{if $k_{i+1}<k_1$},\\
 k_{i+1}-1&\text{if $k_{i+1}>k_1$}.
 \end{cases}
 $$
 \end{proposition}

The proof immediately follows from definitions.
Passing from the sequence $\{k_i\}$ to the sequence $\{r_i\}$ is the ``dynamics'' of our coordinates.

So, we associate with our sequence of partitions $\{\eta_n\}_n$ a tree endowed with additional structures, a translation and a transfer (see Fig.~1). The tree $W$ with these structures will be called the \emph{skeleton} of the permutation tree.

 An ordinary edge of the tree joins a permutation $g \in {\mathfrak S}_n$ with a permutation
$h \in {\mathfrak S}_{n+1}$ if and only if $g$ is the result of removing the $(n+1)$th object from the permutation $h$. In contrast, the translation joins a permutation $g \in {\mathfrak S}_n$
with a permutation $f \in {\mathfrak S}_{n-1}$ if $f$ is the result of  removing the first object from the permutation\footnote{Here, by abuse of notation, we have identified a permutation regarded as an element of the symmetric group ${\mathfrak S}_n$ and  the image of the natural order $1,2, \dots, n$ under this permutation; of course, one must distinguish between these notions, see \cite{V1}.}  $g$.

Thus, we have defined the skeleton, which is a combinatorial scheme of an increasing shift-invariant sequence of finite partitions. Now we want to regard these data as new coordinates of points of the cube.

\begin{definition}[transfer]
Consider the space $T(W)$ of all paths of the permutation tree $W$. The operation $\lambda$ that sends a permutation $g\in {\mathfrak S}_n$ of level $n$ by translation to a permutation $h\in {\mathfrak S}_{n-1}$ of level $n-1$ induces a map $\Lambda:T(W)\rightarrow T(W)$: a path $\{g_1,g_2,g_3, \dots\}$ goes to the path $\{\lambda(g_2),\lambda(g_3), \dots \}$; we call $\Lambda$ the \emph{transfer} on the tree $W$. Since $g\prec h$ implies $\lambda(g)\prec\lambda(h)$   (the two removals commute), this operation defines a mapping on the paths of our tree.
\end{definition}

We will give a transparent formula for this transfer,
see Proposition 2 below.

For a point $x=(x_1,x_2,x_3, \dots) \in I^{\infty} $,
consider the quantities
$t_n(x)=\sharp\{k\leq n: x_k<x_n\}$
and the map
$$
x\mapsto t(x)=\{t_1(x),t_2(x),\dots \},\quad t_n(x)\in \textbf{n},\quad n=1,2 \dots.
$$
It is clear that the sequence $\{t_n(x)\}$ is a path in the skeleton of our tree~$W$.

This is the combinatorial encoding of the infinite-dimensional unit cube $I^{\infty}$ with respect to the partition into Weyl simplices.

\begin{definition}
Consider the compact space
$$\mathfrak{M}=\{\{t_n\}_{n=1}^{\infty}, \quad t_n\in \textbf{n}, \quad n=1,2 \dots\}.
$$
We call $\mathfrak{M}$ the \emph{triangular compact set}.
\end{definition}

We define a map
  $J:I^{\infty}\rightarrow \mathfrak{M}$, $x=\{x_1,x_2, \dots\}\mapsto \{t_1(x),t_2(x) \dots\}$, as follows:
$$
    J(\{x_n\})=\{t_n=t_n(x_1,\dots, x_n)\}_{n=1}^{\infty},
		\quad t_n=\sharp\{i:1\leq i \leq n, \; x_i<x_n\}.
		$$

Also, we define a probability measure  $\mu$ on the compact set $\mathfrak{M}$ as the product
  of the uniform measures on each factor $\textbf{n}$.

\section{The distinguishability of the Weyl encoding}

  Our first observation is as follows.

\begin{theorem}
  The map $J$ is an isomorphism of measure spaces between
  $(I^{\infty}, m^{\infty})$ and $(\mathfrak{M},\mu);$
  this means that the coordinates $\{t_n(x)\}_{n=1}^{\infty}$ determine almost all points $x$, or, in other terms, the product of the partitions~$\eta_n$ is the \emph{identity} partition, or the partition of the infinite-dimensional cube $I^{\infty}$ into  singletons:
$$
\bigvee_{n=1}^{\infty}\eta_n=\epsilon.
$$
\end{theorem}

 In geometric terms, this means that the set of all Weyl simplices of the cube $I^{\infty}$ separates almost all points of the cube. A more paradoxical formulation is as follows: almost every
point  $\{x_1,x_2, \dots\}$ of the cube $I^{\infty}$ with respect to the Lebesgue measure $m^{\infty}$ can be recovered if we know only all inequalities $x_n>x_m$ or $x_n<x_m$ for all $n,m\in \Bbb \N$. Or, in an even more expressive form, \emph{almost every Weyl simplex consists of a single point, and almost every Weyl chamber consists of a single ray}.

\begin{proof}
 Let two sequences $\{x_n\}$ and $\{x'_n\}$ have the same inequalities for all pairs:
 $x_n>x_m \Leftrightarrow x'_n>x'_m$.
 Suppose that $x_1\ne x'_1$; then there exists~$k$
 such that $x_k \in (x_1,x'_1)$ and $x'_k\in (x_1,x'_1)$,
 so the points $x_1, x_k$ satisfy the opposite inequality as compared with $x'_1,x'_k$, a contradiction.
\end{proof}

A more precise form of this claim is as follows.

\begin{lemma}
The limiting partition $\eta=\lim_n \eta_n$  of the infinite-dimensional cube $I^{\infty}$ (the limit of the partitions into open Weyl simplices) coincides $\bmod\, 0$ (with respect to the Lebesgue measure) with the partition into singletons. In other words, the distinguishability problem for the partition into Weyl simplices has a positive solution. Therefore, the map $J$ is an isomorphism of measure spaces. In more detail, there exists a set of full Lebesgue measure in $I^{\infty}$ such that for any two points~$\{\xi_n\}$ and~$\{\xi'_n\}$ of this set there exist indices $i$ and $j$ for which the corresponding coordinates satisfy the opposite inequalities:
$$
\xi_i > \xi_j, \quad \mbox{but}\quad \xi'_i <\xi'_j.
$$
\end{lemma}

We may say that the sequence of partitions into Weyl simplices gives a positive solution to the \emph{distinguishability problem}, the problem of how to separate the points of the infinite-dimensional unit cube via a sequence of finite measurable partitions. The real explanation of this  effect lies in the individual ergodic theorem or in the theorem  on the uniform distribution on the interval of almost all points of the cube. This means that the  Lebesgue measure $m^{\infty}$ can be replaced with any measure having this property.

\smallskip

Now we must describe the translation and transfer defined above in terms of the coordinates $\{t_n\}$, that is, describe how these coordinates change when a point $x=\{x_1,x_2,x_3, \dots\}$ changes into $Sx=\{x_2,x_3, \dots\}$. Namely, we want to express $\{t_1(Sx)=t'_1, t_2(Sx)=t'_2, \dots \}$ in terms of $ \{t_1(x),t_2(x),t_3(x), \dots \}$. Such an expression will give a formula for the transformation
$$
\Lambda=JSJ^{-1}: \mathfrak{M} \rightarrow \mathfrak{M}
$$
 of the triangular compact set $\mathfrak M$.

This map $\Lambda$ is the transfer of the permutation tree and of the triangular compact set regarded as a set of paths of this tree. But we will give a precise formula for $\Lambda$.

Denote $\xi^n=\{\xi_i\}_{i=1}^n$, and let  $d_n(\xi^n)$ be the number of coordinates in the vector $\xi^n$ that are less than $\xi_1$. It is clear that $d_{n+1}$ is either equal to $d_n+1$ if $\xi_{n+1}<\xi_1$, or equal to $d_n$
  if $\xi_{d+1}>\xi_1$. Let us say that
 some positions in the vector $(t_1,t_2, \dots ,t_n)$
 are \emph{special}: the first position $t_1=1$ is special; and if the number of special positions among the first $n$ ones is $d_n(t)$, then $t_{n+1}$ is special if and only if
 $t_{n+1}\leq d_n$, $n=1,2, \dots$.

\begin{proposition}[formula for $\Lambda= JSJ^{-1}$]
We have  $\Lambda(\{t_n\})=\{t'_n\}$
 where
 $$
 t'_n=\begin{cases}
 t_{n+1} &\text{if \  $t_{n+1}$ is special},\\
 t_{n+1}-1&\text{if  \ $t_{n+1}$ is not special.}
 \end{cases}
 $$
\end{proposition}

     The following formula holds:
$$
	t_{n+1}-t'_n =1-(d_{n+1}(t)-d_n(t)).
$$

The proof follows automatically from the previous formulas.

Now consider a formula for the inverse map.  From the formula for $d_n$ we
directly obtain the following.

\begin{theorem}
For almost every trajectory  $\{\xi_n\}_n\in I^{\infty}$ with respect to the measure $m^{\infty}$ on $I^{\infty}$,
$$
\lim_n \frac{d_n}{n}=\xi_1.
$$
 \end{theorem}

So, we can find the first coordinate from the infinite vector $\{t_n\}$. In the same manner we can find the other coordinates $\{x_2,x_3, \dots, x_n\}$ using the transfer.

 This means that we can recover the shift $S$
as a map on $I^{\infty}$ with the help of $\Lambda$.

We establish an isomorphism between two triples
  $$
	(I^{\infty},m^{\infty}, S)\quad \mbox{and}\quad (\mathfrak{S},\mu,\Lambda).
	$$
  In particular, this means that $\Lambda$ is a Bernoulli automorphism of
 $(\mathfrak{S},\mu)$ in the sense of ergodic theory.

By definition, the translation (see Sec.~2.3) associates with every vertex of level~$n$ (for $n>1$) a vertex of level~$n-1$ following the rule according to which the simplex
$\sigma_{x^n}$ of sequences starting from a vector
$x^n=(x_1,x_2, \dots, x_n)$  changes when we remove the first coordinate $x_1$ after the application of the shift $S$, that is, pass to the vector
$Sx^n=(x_2,x_3, \dots, x_n)$. Recall that the permutation
$(k_1,k_2,\dots, k_n)$ corresponding to the simplex $\sigma_{x^n}$
is given by the formula
$$
		k_i=\#\{s \in {\bf n}: x_s<x_i\},\quad i=1,2, \dots, n.
$$

Using this rule, we construct a correspondence with the tree of Weyl simplices, see Fig.~1.

Thus, we have defined a translation which is a~map from the set of permutations of length~$n$ to the set of permutations of length~$n-1$. In the next section, where we compute the transfer for this graph, we use this map and interpret it in a~slightly different way.

\section{A more complicated case: encoding via the RSK correspondence; see \cite{RS,Sn}}

In the previous section, we considered the partition $\eta_n$ of the cube $I^n$ into Weyl simplices and constructed an encoding of the dynamical system $(I^{\infty}, m^{\infty}, S)$ using this kind of partitions. We obtained a combinatorial version of a Bernoulli shift as a transformation on the space of paths of a tree (the permutation tree).

Looking at more complicated examples, assume that  partitions $\alpha_n$ of the cube $I^n$ are coarser than the partitions $\eta_n$ into Weyl simplices (${\alpha_n \prec \eta_n}$), namely, each element of $\alpha_n$ consists of several Weyl simplices. In this case, we can loose the distinguishability property: the limit $\bigvee_n \alpha_n$ can be different from the identity partition. How to refine the distinguishability problem for such sequences of partitions? How to single out the case where $\lim_n \alpha_n=\epsilon$?

As one of the important further examples, we will keep in mind  the situation with the so-called RSK correspondence. A related question was considered in the 1980s in~\cite{KV} and recently in the important papers \cite{RS,Sn}.

In the framework of this paper, the positive answer to the distinguishability problem for the combinatorial encoding is equivalent to the fact that the homomorphism of a Bernoulli shift to the space of infinite Young tableaux with the Plancherel measure and the Sch\"utzenberger shift is an isomorphism. The question of whether this homomorphism is an isomorphism appeared as a result of the paper \cite{KV}, in which a generalization of the RSK correspondence for the infinite case was defined; and this question was recently answered in the papers \cite{RS,Sn}. Here we mention another approach to the isomorphism problem in the spirit of this paper and the article~\cite{V}, as well as the theory of filtrations~\cite{UMN}.

We will use the identification of Weyl simplices with permutations; more exactly, a permutation $g=(i_1,i_2,\dots, i_n)\in {\mathfrak S}_n$ parametrizes the Weyl simplex $\sigma_g$ whose elements $x=(x_1,x_2,\dots, x_n)\in I^n$ have the same ordering of  coordinates as  the permutation $g$.

The main property of the finite RSK correspondence (see \cite{St}) is a set-theoretic isomorphism
between the symmetric group ${\mathfrak S}_n$ and the set of all pairs of Young tableaux with the same Young diagrams:
$$
 {\mathfrak S}_n=\coprod_{\lambda \in {\hat S}_n} T^P_{\lambda}\times T^Q_{\lambda},
$$
where $T_\lambda$ is the set of all Young tableaux of shape~$\lambda$; the indices $P,Q$ mean that a tableau
is either insertion ($P$) or recording ($Q$).

The index of any simplex is a permutation $g\in {\mathfrak S}_n$, but we regard $g$ (as above) as a pair of Young tableaux of the same shape  $\lambda_n$, i.e., $g=(t_g^P,t_g^Q$) (in short, we may say that $g$ is associated with the diagram $\lambda_n$).

Now denote by $\Sigma_{\lambda_n}$ the cylinder set that is the union of all simplices $\sigma_g$ with diagram ${\lambda_n}$.

We define two partitions $\eta_n^{\lambda}$ and $\phi_n^{\lambda}$ of $\Sigma_{\lambda_n}$. The first partition $\eta_n^{\lambda}$ is finite and cylinder, with  elements of the form
$$
C_{t_Q}=\bigcup_{t_P}\sigma_{t_P,t_Q},\quad  t_P,t_Q\vdash \lambda.
$$

The second partition $\phi_n^{\lambda}$ is
cofinite (i.e., all its elements are finite sets), each its element  is a finite set of points  $\{x_n\}\in I^{\infty}$ with equal coordinates $x_i$ for $i>n$,
and the vector $(x_1,x_2,\dots, x_n)$ runs over all normalized permutations with given tableau $t_P$:
$$
D_{t_P}=\bigcup_{g=(t_P,t_Q)}\Big\{x=(x_1,x_2,\dots, x_n): x_i=\frac{g(i)}{n},\; i=1,2, \dots, n\Big\}.
$$

\begin{lemma}
The partitions  $\eta_n^{\lambda}$  and $\phi_n^{\lambda}$  of the set $\Sigma_{\lambda_n}$ are mutually independent complements: every element of the first partition and every element of the second partition intersect in exactly one point, and the independence properties for the conditional measures also hold. The partition $\eta_n^{\lambda}$ increases with $n$, and the partition $\phi_n^{\lambda}$ decreases with $n$.
\end{lemma}

So, the  elements of the partition $\eta_n^{\lambda}$ (respectively, $\phi_n^{\lambda}$) are parametrized by the set of  permutations with given $P$-tableau (respectively, with given $Q$-tableau). Recall that in the previous section, an element of the partition~$\eta_n$ was parametrized by one permutation.

Now we can consider the partition $\theta_n$ (respectively, $\theta_n^{\bot}$) of the cube~$I^{\infty}$ whose restriction to each set $\Sigma_{\lambda_n}$ is $\eta_n^{\lambda}$ (respectively, $\phi_n^{\lambda}$).  It is clear that the first
 partition is finite and measurable (the number of elements is equal to the number of Young tableaux with  given number of cells), and the second one is cofinite (the number of points in one element of the partition is finite and does not depend on the element). Obviously,  the partition $\theta_n$ is coarser than $\eta_n$ (i.e., $\theta_n \prec \eta_n$), and the partition $\theta_n^{\bot}$ is finer than the orbit partition of the symmetric group ${\mathfrak S}_n$ for all $n$. Moreover,
$\theta_n \prec \theta_{n+1}$, $n=1,2, \dots$, and
$\theta_n^{\bot}\succ \theta_{n+1}^{\bot}$, $n=1,2, \dots$.

It is useful to give an interpretation of these partitions
in terms of the Knuth equivalence and dual Knuth equivalence, see the book~\cite{St} and Fomin's addendum to the Russian translation of this book.

\medskip\noindent
{\bf Question.} What is the limit of the increasing sequence of partitions $\theta_n$?
\smallskip

The following theorem is equivalent to the remarkable theorem  by D.~Ro\-mik and
P.~Sniady \cite{RS,Sn} in which they proved that the homomorphism defined in \cite{KV} is a true isomorphism. This homomorphism in the case of the Plancherel measure was defined as a map from $I^{\infty}$ to the space of infinite standard Young tableaux. It sends the Lebesgue measure to the Plancherel measure. The question of whether this homomorphism is indeed an isomorphism remained open for many years and was solved in \cite{RS,Sn}. We give a statement of the theorem in terms related to our approach.

\begin{theorem}
 The limit of the sequence of partitions $\theta_n$ on the space $I^{\infty}$ is the identity partition:
$$
\lim_n \theta_n=\epsilon.
$$
Thus, the distinguishability problem (see above) has a positive solution.
\end{theorem}

\begin{corollary}
The limit of the decreasing sequence of partitions $\theta_n^{\bot}$ is the trivial partition:
$$\lim_n \theta_n^{\bot}=\nu.$$
\end{corollary}

Note that, in general, the fact that the limit of a decreasing sequence of partitions $\{\alpha_n\}$ is trivial is only necessary, but not sufficient for  the limit of an increasing sequence $\{\alpha_n^{\bot}\}$ of partitions that are independent complements to $\alpha_n$ to be the identity partition.

In terms of the dual Knuth equivalence relation on the cube $I^{\infty}$, the conclusion of the theorem means that
it is the identity equivalence relation $\bmod\, 0$.

At the same time, we have the following corollary of the theorem above.

\begin{proposition} The direct Knuth equivalence relation {\rm(}see {\rm\cite{St})} on the space $I^{\infty}$ is ergodic; this means that there is no nonconstant measurable function that is constant on the classes of the dual Knuth equivalence.
\end{proposition}

Our numerical simulations have shown a very slow convergence in this situation.

The approach of the papers \cite{RS,Sn} is based on the analysis of the limit shape of Young diagrams and the so-called Sch\"utzenbeger transformation (``jeu de taquin''). More precisely, each infinite Young tableau is a path in the Young graph; the Sch\"utzenbeger transformation is a shift, or, in our terminology, a transfer of the path (see the next section). It is  natural to choose the maximal strictly monotone subset of the path, which is called the ``nerve'' of the tableau. The main observation of the authors of \cite{RS,Sn} is that this nerve, regarded as a piecewise linear line on the lattice ${\Bbb Z}^2_+$, has a limit at infinity, which can be identified  (after a normalization) with a point of the limit shape of the normalized Young diagram. This is the key step of constructing the inverse isomorphism in \cite{RS,Sn}.

\smallskip
Our approach is different and based on another limit shape theorem. We will give the details elsewhere. Namely, our proof uses a detailed analysis of the behavior of $P$-tableaux. These $P$-tableaux are not semistandard, because their entries are reals but not integers. In contrast to the case of $Q$-tableaux, the $P$-tableaux have no strong limit (as $n$ tends to infinity), and their weak limit is the zero tableau. Nevertheless, after a proper normalization, they have a limit shape; and its existence can be extracted from old results on the  limit shape of standard Young tableaux (not diagrams). The inversion formula from \cite{RS} is based on the  behavior of $Q$-tableaux. Our approach to the inversion formula is based on the stabilization of $P$-tableaux for a given $Q$-tableau after the normalization. Roughly speaking, there is an isomorphism between the realizations of the Bernoulli scheme ($I^{\infty}$), the infinite Young tableaux ($Q$-tableaux with the Plancherel measure), and the
normalized $P$-tableaux.

\bigskip

\section{Definition of transfer for a graded graph}

Let us define the notion of ``transfer,''  an operation on the path space of a graded graph (or multigraph). The main idea is to apply the theory of graded graphs to ergodic problems under consideration.

Consider an arbitrary Bratteli diagram, i.e., an $\Bbb N$-graded locally finite graph (or even a multigraph)~$\Gamma$. An infinite tree is an example of such a graph. A path in~$\Gamma$ is an infinite maximal sequence of edges (not vertices!) in which the beginning of each edge coincides with the end of the previous edge. Denote the space of all paths by  $T(\Gamma)$; this is a Cantor-like compactum in the inverse limit topology.

\begin{definition}
A transformation
$$
\Lambda: T(\Gamma)\rightarrow T(\Gamma),\qquad \Lambda(\{t_n\}_n)=(\{u_n\}_n),
$$
is called a \emph{transfer of general type} if it is a continuous map in the compact topology of the path space, and it is defined by the following successive local rules:
   for every $n$, the edge $u_n$ of the image depends on a fragment of the argument path and the previous edge of the path image: $u_n= f_n(\{t_i\}_1^{n+1}, u_{n-1})$, where $f_n$ are arbitrary combinatorial functions, $n=1,2,\dots$.
   This definition describes a very general class of transfers. For trees regarded as graded graphs, it includes the definition from Sec.~{\rm1}.
\end{definition}

For stationary graphs, in which the sets of vertices of every level (except the first one) are isomorphic and these isomorphisms are fixed, the translation rule depends on nothing: an edge connecting vertices~$a$ and $b$ of levels~$n+1$ and~$n+2$ goes to the edge connecting the vertices~$a'$ and~$b'$ of levels~$n$ and~$n+1$ identified with the vertices~$a$ and~$b$, respectively. In this case, the transfer is an ordinary shift. We say that a transfer is Markov if the fragment $\{t_1,\dots, t_{n+1}\}$ above reduces to the last pair of  edges $\{t_n,t_{n+1}\}$.

So, a transfer gives a nonstationary generalization of stationary models of dynamics.

For graphs of a special type, the definition of transfer can be formulated in a very simple way. Assume that a graded graph $\Gamma$ has the following property: for every $2$-interval $[v,w]$ of the graph, with $\operatorname{grad}(v)=n$, $\operatorname{grad}(w)=n+2$, there is an involution $\phi(v,w)\equiv \phi$ on the set of vertices of the interval that have level $n+1$ (we omit the index of $\phi$ when it is clear from the context). Then if $v=(v_1,v_2,\dots )$ is a path, then the transfer is given by the following formula:
$$
\Lambda(v)=(v_1,\phi(v_2),\phi(v_3),\dots).
$$
In this case, the natural transfer is Markov, see Fig.~2 from which our definition should be clear.

For all distributive lattices, the transfer is well defined, because all $2$\nobreakdash-intervals in this case have one or two intermediate vertices. In particular, this definition is valid for the Young graph.

\begin{definition}
A graded graph with a ``transfer'' operation is called a quasi-stationary graph. The path space of a~quasi-stationary graph, regarded as a topological Markov compactum, is called a quasi-stationary Markov compactum.
\end{definition}

Thus, we have described a new type of realizations of automorphisms and endomorphisms with infinite entropy as transfers on quasi-stationary Markov compacta, i.e.,
in spaces that are locally finite.

According to our definition, a transfer is a shift of sequences of edges, and not of sequences of vertices as in the stationary case. Thus, this notion opens new possibilities for realizations of transformations.

\begin{figure}
  \includegraphics[scale=0.5]{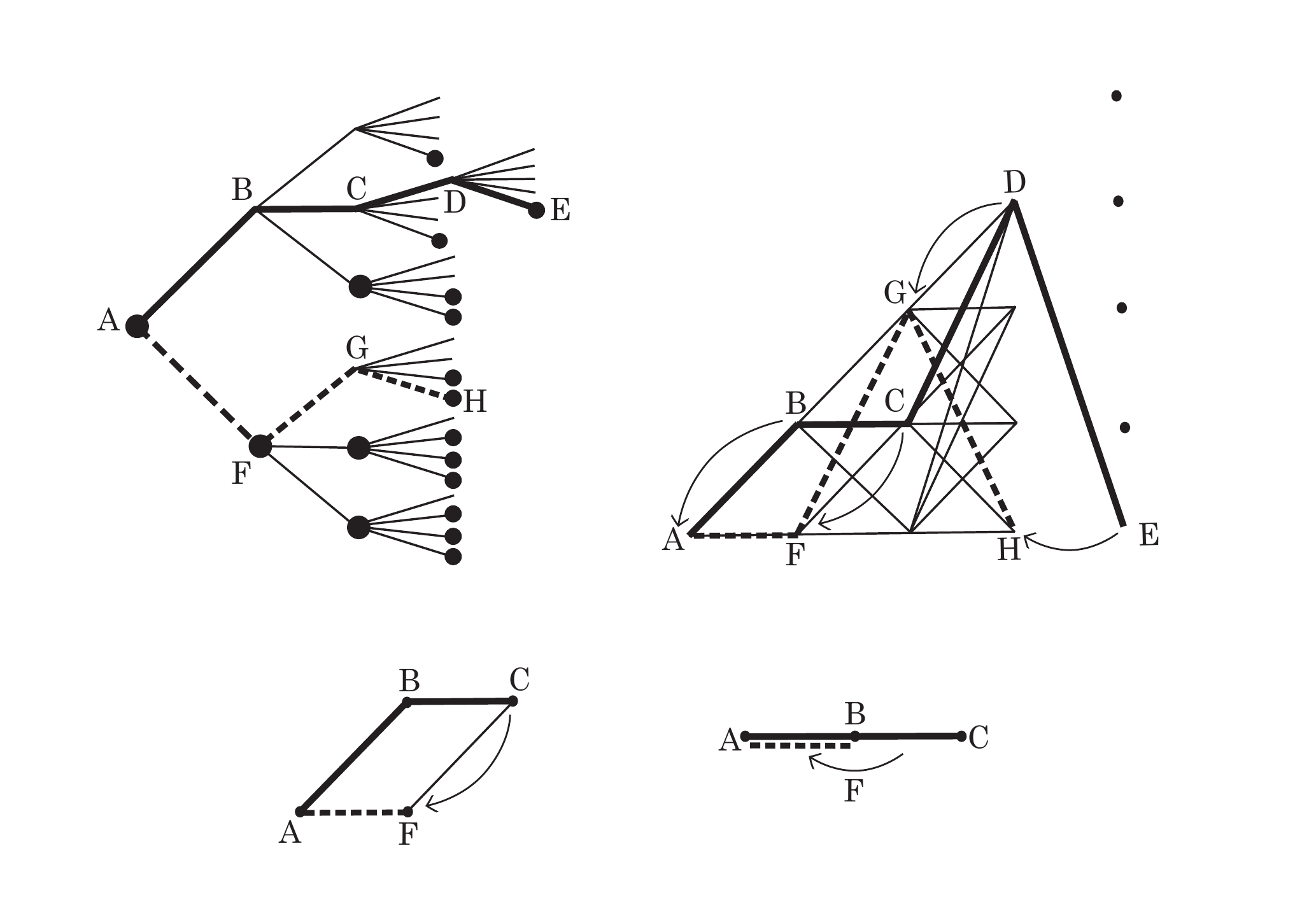}
  \caption{The dashed path is the transfer of the bold one.}
  \label{fig:transfer}
\end{figure}


\begin{proposition}
For a graph that is the Hasse diagram of the distributive lattice of finite ideals of a locally finite countable partially ordered set with a minimal element, there is a distinguished Markov transfer. In the case of the lattice ${\Bbb Z}_2$ as a poset and the Young graph as a graded graph, the transfer coincides with the well-known Sch\"utzenberger transformation {\rm(}``jeu de taquin,'' see {\rm\cite{St})}, so our definition is a generalization of ``jeu de taquin'' for distributive lattices.
\end{proposition}

The proof follows from a detailed analysis of the definition of transfer.
If  a transfer is defined on the path space of a graded graph, then this space should be regarded as a~nonstationary (or {\it quasi-stationary}) Markov chain, meaning that  the transfer is an analog of the shift. If we have a central measure on the path space that is invariant under the transfer, then we obtain a~quasi-stationary Markov chain with an invariant measure.  Hence the theory of transfer becomes part of ergodic theory, as a nonconventional realization of measure-preserving transformations. We will return to all these facts elsewhere.

\bigskip
The author is grateful to P.~Nikitin for preparing the figures and for the design of the manuscript.

\end{document}